\documentclass{amsart}
\newtheorem{thm}{Theorem}
\newtheorem{cor}[thm]{Corollary}
\newtheorem{lem}[thm]{Lemma}
\newtheorem{prop}[thm]{Proposition}
\theoremstyle{definition}

\newtheorem{exmp}{Example}
\newtheorem{remark}{Remark}

\newcommand{\BEX}{\begin{exmp}}
\newcommand{\EEX}{\end{exmp}}

\allowdisplaybreaks[4]
  \def\Rset{\mathbb{R}}

\newcommand{\nwc}{\newcommand}
\nwc{\COM}[1]{}

\newcommand{\neprop}[1]{{Proposition \ref{#1}}}
\newcommand{\netheo}[1]{{Proposition \ref{#1}}}

\newcommand{\kb}[1]{\boldsymbol{#1}}
\newcommand{\vk}[1]{\kb{#1}}

\def\a{\vk{a}}

\newcommand{\ve}{\varepsilon}
\newcommand{\abs}[1]{\lvert #1 \rvert}

\newcommand{\pk}[1]{\mbox{\rm$\vk{P}$} \{#1\} }

\newcommand{\pb}[1]{\mbox{\rm$\vk{P}$}\Bigl \{#1 \Bigr \}}

\def\R{\Rset}

\newcommand{\inr}{\in \R}

\newcommand{\limit}[1]{\lim_{#1 \to   \infty}}
\newcommand{\todis}{\stackrel{d}{\to}}

\newcommand{\equaldis}{\stackrel{d}{=}}

\newcommand{\BQN}{\begin{eqnarray}}
\newcommand{\EQN}{\end{eqnarray}}
\newcommand{\BQNY}{\begin{eqnarray*}}
\newcommand{\EQNY}{\end{eqnarray*}}
\newcommand{\BS}{\begin{prop}}
\newcommand{\ES}{\end{prop}}
\newcommand{\BRM}{\begin{remark}}
\newcommand{\ERM}{\end{remark}}

\newcommand{\BL}{\begin{lem}}
\newcommand{\EL}{\end{lem}}
\newcommand{\BT}{\BS}
\newcommand{\ET}{\ES}
\newcommand{\BK}{\begin{cor}}
\newcommand{\EK}{\end{cor}}

\newcommand{\QED}{\hfill $\Box$}
\newcommand{\IF}{\infty}

\newcommand{\proofprop}[1]{\textsc{Proof of Proposition} \ref{#1}}

\def\Ay{a_\rho}
\def\Ayx{x_*}
\def\Orho{\rho_*}

\def\alphAR{\vk{\alpha}_{a, \rho}}
\def\alphAR{\Ay}

\def\uY{\Ay x}

\def\paxD{p_{a,\delta, \eta; x}}
\def\paxDD{p_{a,\delta, \eta, \rho; x}}

\begin{document}
\title[Exact Tail Asymptotics in Bivariate Scale Mixture Models]{Exact Tail Asymptotics in Bivariate Scale Mixture Models}

\author{Enkelejd Hashorva}
\address{Department of Mathematical  Statistics and Actuarial Science\\
University of Bern, Sidlerstrasse 5\\
CH-3012 Bern, Switzerland\\
enkelejd.hashorva@Stat.Unibe.ch }


\date{\today{}}
\begin{abstract}
Let $(X,Y)=(RU_1, R U_2)$ be a given bivariate scale mixture random vector, with  $R>0$ being independent of the bivariate random vector $(U_1,U_2)$.
 In this paper we derive exact asymptotic expansions of the joint survivor  probability of $(X,Y)$
  assuming that $R$ has distribution function in the Gumbel max-domain of attraction and $(U_1,U_2)$ has a specific tail behaviour  around some absorbing point.   We apply our results to investigate the asymptotic behaviour of joint conditional excess distribution and
  the asymptotic independence for two models of bivariate scale mixture distributions. Furthermore for our models we derive an expression of the residual dependence index $\eta$.
\end{abstract}


 \maketitle


\section{Introduction}
Let $(X,Y)$ be a bivariate random vector with stochastic representation
\BQN \label{eq:S}
(X, Y) \equaldis (RU_1,R U_2),
\EQN
where $R>0$ is independent of the bivariate random vector $(U_1,U_2)$ ($\equaldis$ stands for equality of distribution functions). 
The random vector $(X,Y)$ has  a scale mixture distribution; a canonical example of such $(X,Y)$ is a bivariate spherical random vector with rotational invariant distribution function with $(U_1,U_2)$ uniformly distributed on the unit circle of $\R^2$.
In this model (see Cambanis et al.\ (1981)) the dependence between $U_1$ and $U_2$ is a functional one, namely $U_1^2+U_2^2=1$ almost surely, and
\BQN \label{eq:sph}
 (U_1, U_2  ) &\equaldis &(I_1 W, I_2\sqrt{1- W^2} ),
 \EQN
with $W\in (0,1), I_1,I_2 \in \{-1, 1\}$ almost surely, $W^2$  beta distributed with parameters $1/2,1/2$ and $\pk{I_1=1}= \pk{I_2=1}=1/2$.
Furthermore,  $I_1,I_2,R, W$ are mutually independent.

If $R$ is such that $R^2$ is Chi-squared distributed with two degrees of freedom, then $X$ and $Y$ are independent Gaussian random variables with mean zero and variance 1.\\
  Our main interest in this paper is the tail asymptotics of the joint survivor function of $(X,Y)$.
For Gaussian random vectors the asymptotics of the joint survivor  probability is well-known,
see e.g., Berman (1962), Dai and Mukherjea (2001), Hashorva (2005), or Lu and Li (2009).
Results for elliptical and Dirichlet random vectors are obtained in Hashorva (2007, 2008, 2009c) and  Manner and Segers (2009).
Note that the elliptical model is derived by extending \eqref{eq:sph} to
\BQN \label{eq:ell}
 (U_1, U_2  ) &\equaldis &(I_1 W, I_1\rho W+  I_2\Orho \sqrt{1- W^2} ), \quad \rho \in (-1,1), \quad \Orho:=\sqrt{1- \rho^2}, \EQN
 where the additional parameter $\rho$ corresponds to correlation coefficient of $X$ and $Y$ if $R^2$ is Chi-squared distributed.\\
Hashorva (2007) generalises the known asymptotic results for Gaussian random vectors to the more general class of elliptical ones by
exploiting the fact that the asymptotics of the joint survivor probability is primarily determined by the asymptotic properties of the survivor function  $\overline{F}:= 1- F$ of the associated random radius $R$.\\ 
Specifically, in the aforementioned paper the principal asymptotic assumption is that $F$ is in the Gumbel max-domain of attraction, which means that for some positive scaling function $w$
 \BQN
\label{eq:gumbel:w}
\limit{x} \frac{\overline{F}(x+t/w(x))}{\overline{F}(x)}&=& \exp(-t), \quad \forall t\inr.
\EQN
As shown in Hashorva (2007)  condition \eqref{eq:gumbel:w} is crucial when $(X,Y)$ is
an elliptical random vector with stochastic representation \eqref{eq:S}.
 More specifically, by the aforementioned paper for any $a\in (\rho, 1]$
\BQN\label{eq:2}
\pk{X> x, Y> ax}&\sim&   \frac{\Ay^2\Orho^3}{2 \pi (1-a\rho)(a- \rho)}   \frac{1}{v(\Ay x) }
\overline{F} (\Ay  x), 
\EQN
where
\def\vx{v(x)}
\BQN \label{Ay:vx}
\Ay&:=& \Orho^{-1}\sqrt{1- 2a\rho+ a^2}>1, \quad \vx:=x w(x), \quad x\inr.
\EQN
Throughout this paper $f(x)\sim g(x)$ means $\limit{x} f(x)/g(x)=1$, and \eqref{eq:gumbel:w} is abbreviated by $F \in GMDA(w)$.

If $W^2$ is beta distributed with positive parameters $\alpha,\beta$ (its distribution function is denoted by
$beta(\alpha,\beta)$), then $(X,Y)$ is a generalised Dirichlet random vector. Hashorva (2009c) extends \eqref{eq:2} for the class of Dirichlet random vectors.
As indicated in Hashorva (2009b) for certain asymptotic problems the distributional properties of $(U_1,U_2)$ do not need to be explicitly known. A natural question that arises concerning the asymptotics of the joint survival probability of $X,Y$ is that if  $F\in GMDA(w)$ what models for the dependence between $(U_1,U_2)$ would lead us to asymptotic results similar to \eqref{eq:2}?

In this paper we answer the above question for two specific models:  The first one is refereed to as the unconstrained dependence model, or simply Model A. In that model we assume that $U_1\in (0,1]$ almost surely, and further impose an asymptotic assumption on the behaviour of the $(U_1,U_2)$ around some absorbing point $(1,a)$ (see \eqref{theo1:1:1} below).

The second model (or simply Model B) motivated by \eqref{eq:ell} is referred to as the functional dependence model. More specifically we assume
the stochastic representation
\BQN \label{eq:wZ}
 (U_1, U_2  ) &\equaldis &(I_1 W, \rho I_1 W + I_2z^*(W)), \quad \rho \in (-1,1),
 \EQN
with $z^*$ some positive measurable function, $W\in (0,1), I_1,I_2\in \{-1,1\}$ almost surely,
and $I_1,I_2, W$ are mutually independent.

We present three applications of our results: The first one establishes an asymptotic approximation of the joint conditional excess distribution. In the second application we discuss the Gumbel max-domain of attraction of bivariate distributions related to our Model B.
In our last application we derive an explicit expression of the residual tail dependence index $\eta$ for bivariate scale mixture random vectors extending a recent result of Hashorva (2010) for elliptical random vector.

Organisation of the paper: In the next section  we state our first result dealing with some general scale mixture bivariate random vectors which fall under
Model A. We introduce in Section 3 some constrains on the dependence function
of $(U_1,U_2)$ via \eqref{eq:wZ}, and then investigate the tail asymptotics of interest for Model B showing a
 generalisation of \eqref{eq:2} in \neprop{prop:AAA}. Three  applications of our results are given in Section 5.
 Proofs of all the results are relegated to Section 6.

\def\wpax{\widetilde p_{a}(s)}
\section{Tail Asymptotics Under Unconstrained  Dependence}
Consider a bivariate scale mixture random vector $(X,Y)=(RU_1, R U_2)$, where $R$ has distribution function $F$ (denote this $R\simeq F$)
satisfying \eqref{eq:gumbel:w} with some positive scaling function $w$. We assume throughout this paper that $F$ has an infinite upper endpoint.  Hence by \eqref{eq:gumbel:w}   (see e.g., Resnick (2008))
\BQN\label{eq:uv}
\vx:=x w(x) &\to &\IF, \quad x\to \IF.
\EQN
Given a constant $a\in(0,1]$ we investigate the asymptotics of
$$\paxD:=\pk{X>x[1 + \delta/\vx], Y> a x[1+ \eta/\vx]}, \quad x\to \IF,$$
 for any $\delta,\eta \in [0,\IF)$. The reason for dealing with the asymptotics of
 $\paxD$ is our interest concerning the approximation and estimation of the joint conditional excess distribution, see the first application in Section 4.\\
Throughout in the sequel we assume that $U_1$ is a bounded random variable. Without loss of generality we consider only the case $U_1$ has distribution function with upper endpoint equal 1. This implies that $\paxD\le \overline{F}(x)$ for any $x$ positive.
For both Model A and B  we  show below that this upper bound is too crude; roughly speaking we have the asymptotic behaviour
 $$ \paxD\sim \psi(x) \overline{F}(x),$$
with $\psi$ some positive function decaying polynomially fast to 0 as $x\to \IF$.

In addition to the Gumbel max-domain of attraction assumption, we impose next a certain asymptotic behaviour of  $(U_1,U_2)$ around $(1,a)$, namely
\BQN\label{theo1:1:1}
\limit{x} \frac{\pk{U_1 > 1- (s-\delta)/x, U_2 > a(1- (s- \eta)/x) }}
{\pk{U_a > 1- 1/x }}  &=& \xi_a(s,\delta,\eta), \quad \forall \delta, \eta \in [0,\IF), \forall s\in (0,\IF),
\EQN
with $\xi_a$ a positive measurable function and $U_a:=\min(U_1,U_2/a)$. If $\delta\ge \eta$, then
$$ \xi_a(s,\delta, \eta)= \xi_a(s, \delta- \eta, 0), \quad \xi_a(s,\eta, \delta)= \xi_a(s, 0, \delta- \eta), \quad \forall s\in (0,\IF).$$
Further, for $\delta=\eta=0$ 
\BQN \label{eq:paxg}
\pk{U_a > 1- s}=s^\gamma L_a(s), \quad  \forall s>0
\EQN
holds for some $\gamma \in [0, \IF)$, with $L_a$ a positive measurable function such that 
$\lim_{s\downarrow  0} L_a(s)/L_a(ts)=1,\forall t>0$, i.e., $L_a$ is a slowly varying function;
see Bingham et al.\ (1987), Embrechts et al.\ (1997), Falk et al.\ (2004), De Haan and Ferreira (2006),
Jessen and Mikosch (2006), Resnick (2008), or Omey and Segers (2009)  for more details on regularly varying functions.

 Next, we formulate our first result.

\BT \label{theo1} Let $(X,Y)=(R U_1, R U_2)$ be a bivariate scale mixture random vector with
 $R\simeq F$  a positive random variable being independent of $(U_1, U_2) $.  Suppose that $F$ has an infinite upper endpoint satisfying \eqref{eq:gumbel:w} with some positive scaling function $w$, and $U_1\in (0,1]$ has distribution function with upper endpoint 1.
 If  $a\in (0, 1]$ is such that \eqref{theo1:1:1} holds, then for any $\delta,\eta \in [0,\IF)$ we have
\BQN\label{theo1:2}
\paxD
&\sim & J_{\delta, \eta} L_a(1/\vx) \frac{\overline{F}(x)}{(\vx)^{-\gamma } }  ,
\EQN
with $L_a$ satisfying \eqref{eq:paxg} and
$$ J_{\delta, \eta}  := \int_\delta^\IF \xi_a(s,\delta, \eta) \exp(-s)\, ds \in (0,\IF).$$
\ET

{\bf Remarks:}  (a) In view of Lemma 6.1 in Hashorva (2009b)
for any $\lambda\in (1,\IF),  c\inr$ and $F$ as in \netheo{theo1} we have
\BQN\label{eq:resn}
\limit{x} \frac{(\vx) ^c \overline{F}(\lambda x)}{ \overline{F}(x)}&=&0.
\EQN
In fact \eqref{eq:resn} follows directly from Proposition 1.1 in  Davis and Resnick (1988), see also
Embrechts et al.\ (1997) p.\ 586, and {\bf A1.} in Hashorva (2009c).

Further we have the self-neglecting property of $w$, i.e.,
 \BQN\label{eq:wb}
\frac{w(x+ t/w(x))}{w(x)}&\sim &1
\EQN
holds locally uniformly for $t\inr$.  Refer to Galambos (1987), Reiss (1989), Embrechts et al.\ (1997), Falk et al.\ (2004),
de Haan and Ferreira (2006), or Resnick (2008) for  details on the Gumbel max-domain of attraction.\\
(b) Under the assumptions of \netheo{theo1} it follows that
\BQN
 \paxD &\sim& \overline H(x),
 \EQN
 where $H$ is the  distribution function of $R W_{a}$ with $W_{a}$ a positive random variable
 independent of $R$ such that
 $$\pk{W_a> 1- s}=J_{\delta,\eta} L_a(s)s^\gamma, \quad s\in (0,1).$$
 See also Example 1 below.  Furthermore \eqref{theo1:2} holds locally uniformly in $\eta,\delta$.

(c) Since for $\delta=\eta=0$ \eqref{eq:paxg} holds, then $J_{0,0}= \Gamma(\gamma+1).$
By the monotonicity of $J_{\delta,\eta}$ in $\delta,\eta$ we obtain
$$ J_{\delta,\eta} \le \Gamma( \gamma+1), \quad \forall \delta,\eta>0.$$

(d) By \eqref{eq:gumbel:w}, \eqref{theo1:2} and   \eqref{eq:wb} it follows that
$$ J_{\delta, \eta} = \exp(- \delta) J_{0, \eta- \delta}, \quad \forall \eta\in [0,\delta],$$
which follows also directly by the definition on $J_{\delta, \eta}$ and \eqref{theo1:1:1}.

We present next three illustrating examples.

{\bf Example 1.} Let $U_1, U_2$ be two random variables taking values in $[0,1]$ such that $U_2 \ge U_1$ almost surely.
Suppose that $\pk{U_1> 1 - s}= s^\gamma L(s), s\in (0,1)$ with $\gamma \in [0,\IF)$ and $L$ a
slowly varying function at 0.  Since for any $x>1$
$$ \pk{U_1> 1- 1/x, U_2 > 1- 1/x} =\pk{U_1 > 1- 1/x}$$
if $R$ is independent of $(U_1,U_2)$ satisfying the assumptions of \netheo{theo1}
we obtain
\BQN \label{eq:berm}
\pk{ RU_1 > x}=\pk{ RU_1 > x, R U_2 >  x }
&\sim&\Gamma(\gamma+1) L(1/\vx) \frac{\overline{F}(x)}{(\vx)^{\gamma}} .
\EQN
We note that for $U_1^2\simeq  beta(\alpha,\beta)$ the asymptotics in \eqref{eq:berm} is shown in Berman (1983), see also Berman (1992).
For the more general case that $U_1$ has a regularly varying survivor function see Theorem 3.1 in Hashorva et al.\ (2009). \\

{\bf Example 2.} (Linear Combinations)
Let $S_i\simeq  G_i, i=1,2$ be two independent random variables with values in $[0,1]$ such that
\BQN\label{eq:Gi}
\limit{x} \frac{\overline{G}_i(1-s/x) }{\overline{G}_i(1-1/x)}&=& \gamma_i, \quad \forall s>0, i=1,2,
\EQN
with $\gamma_i \in [0,\IF)$. Let $\lambda_1, \lambda_2 \in (0,1)$ be given constants with
$ \lambda_1 \ge \lambda _2 $, and set
$$ U_i:= \lambda_i S_1+ \overline{\lambda_i} S_2,  \quad  \overline{\lambda}_i:= 1- \lambda_i, \quad i=1,2.$$
By the definition both $U_1,U_2$ have upper endpoint $1$. For any $\delta,\eta \ge 0$ we have (the proof is postponed to the last section)
\BQN \label{eq:apend:1:00}
\pk{U_1 > 1- (s- \delta)/x, U_2 > 1- (s- \eta)/x} &=& \tilde \xi (s, \delta, \eta) \prod_{i=1}^2 \overline{G}_i(1- 1/x) , \quad s>\max(\delta,\eta),
\EQN
with
$$ \tilde \xi (s,\delta, \eta):=\gamma_2\int_0^\IF \Bigl(\max\Bigl(0,  \min(   [s - \delta - \overline{\lambda_1}z]/\lambda_1,
 [s - \eta- \overline{\lambda_2}z]/\lambda_2)\Bigr) \Bigr)^{\gamma_1} z^{\gamma_2-1} \, dz .$$
Note that for $s\in (0, \max(\delta,\eta)]$ \eqref{eq:apend:1:00} holds with $\tilde  \xi(s,\delta,\eta)=0$, and
when $\delta=\eta=0$
\BQN \label{eq:apend:1}
\tilde \xi ( s,0,0)  &\sim&
C_{\gamma_1,\gamma_2, \lambda_1, \lambda_2}s^{\gamma_1+ \gamma_2}, 
\EQN
with $C_{\gamma_1,\gamma_2, \lambda_1, \lambda_2}\in (0,\IF)$
given by
$$C_{\gamma_1,\gamma_2, \lambda_1, \lambda_2}:= \lambda_1^{- \gamma_1}\int_0^1 [1- \overline{\lambda_1} t]^{\gamma_1} t^{\gamma_2 -1}\, dt
 + \lambda_2^{- \gamma_1}\int_1^{ 1/\overline{\lambda}_2}  [1- \overline{\lambda_2}t]^{\gamma_1} t^{\gamma_2 -1}\, dt   .$$
Consequently, \eqref{theo1:1:1} holds with
$$ \xi_1(s, \delta, \eta):= C_{\gamma_1,\gamma_2, \lambda_1, \lambda_2} \tilde \xi(s, \delta, \eta) s^{- \gamma_1-\gamma_2}\vk{1}_{(s> \max(\delta, \eta))},
 \quad s>0,$$
where $\vk{1}_{( )}$ is the indicator function. Thus with $R\simeq F$ such that $F\in GMDA(w)$ the result of \netheo{theo1} holds.

{ \bf Example 3.} (Farlie-Gumbel-Morgenstern Dependence) Let $ U_i\simeq G_i, i=1,2$ be two random variables with values in $[0,1]$. Suppose that for some $K\in [0,1)$
$$ \pk{U_1> x, U_2> y} = \overline{G}_1(x)  \overline{G}_2(y)[1+ K \overline{G}_1(x)\overline{G}_2(y)], \quad \forall x,y\in [0,1].$$
The bivariate random vector $(U_1,U_2)$ possesses thus the Farlie-Gumbel-Morgenstern distribution (see for more details Hashorva and H\"usler (1999)). If \eqref{eq:Gi} holds, then for any $\delta,\eta,s\in (0, \IF)$ we obtain
$$ \pk{U_1> 1- (s-\delta)/x, U_2> 1 - (s- \eta)/x} \sim (s- \delta)_+^{\gamma_1} (s- \eta)_+^{\gamma_2}  \prod_{i=1}^2 \overline{G}_i (1- 1/x), $$
with $(x )_+:=\max(x,0), x\inr$.
Consequently, if the positive random variable $R\simeq F$ is independent of $(U_1,U_2)$ and $F\in GMDA(w)$, then locally uniformly in $\delta,\eta$
\BQN\label{FGM}
\paxD &\sim&
\int_0^\IF(t- \delta)_+^{\gamma_1} (t-\eta)_+^{\gamma_2}\exp(-t)\, dt \prod_{i=1}^2\overline{G}_i(1- 1/\vx)\overline{F}(x).
\EQN
For any $a\in (0,1)$
we observe another asymptotic behaviour, namely
$$ \pk{U_1> 1- (s-\delta)/x, U_2> a( 1 - (s- \eta)/x)} \sim (s- \delta)_+^{\gamma_1} \overline{G}_1(1- 1/x)\overline{G}_2(a). $$
Consequently,
\BQNY
\paxD &\sim&
\overline{G}_2(a) \Gamma(\gamma_1+1) \exp( - \delta) \overline{G}_1(1- 1/\vx)\overline{F}(x).
\EQNY

\section{Tail Asymptotics For Functional Dependence}
In this section we deal with bivariate scale mixture random vectors assuming that the dependence between the components is
determined by some deterministic function. Explicitly, let $(X,Y)$ be a  bivariate random vector with stochastic
representation
\BQN \label{eq:wzW}
(X,Y)&\equaldis &(R I_1 W, \rho RI_1 W+ R I_2 z^*(W)), \quad \rho \in (-1,1),
\EQN
with $(I_1, I_2), R>0, W\in (0,1)$ mutually independent, and  $z^*:[0,1] \to [0,1]$ a positive measurable function.\\
We assume that the distribution function $F$ of $R$ has an infinite upper point, and that of $W$ has upper endpoint equal 1. In the sequel $I_1,I_2$ take values in
$\{-1,1\}$ with $\pk{I_1=I_2=1}\in (0,1].$ We allow $I_1$ and $I_2$ to be independent.\\
  The random vector $(X,Y)$ is a scale mixture random vector for which the dependence of the components is being determined
  by $\rho, z^*$ and the random variables $R,W,I_i,i=1,2$.
We refer  to the implied dependence of the components as the functional dependence.  Note in passing that if
$$\quad W^2 \simeq beta(1/2,1/2), \quad z^*(x)= \sqrt{1- x^2}, \quad x\in [0,1],$$
and $I_1,I_2$ are independent assuming values $-1,1$ with equal probability,  then $(X,Y)$ is an elliptical random vector.

Generally speaking,  under the setup of \eqref{eq:wzW} it turns out that the local asymptotics of the density function of $W $
is important. More precisely, we are able to provide an asymptotic expansion of
$$p_{a;x}:=p_{a,0,0;x}= \pk{X> x, Y> ax}, \quad a\in (0,1], x>0$$
requiring further that
\BQN \label{eq:cc}
\pk{ W- 1/\Ay \in (K_1 u, K_2 u)}&\sim& L_{K_1,K_2}(u) u^{\gamma_a},  \quad \gamma_a\in [0,\IF)
\EQN
holds for all $u>0$ small with  $K_1< K_2, K_1,K_2\inr$ some given constants such that
$L_{K_1,K_2}(u), u>0$ is a locally bounded  slowly varying function at $0$. Additionally we need to impose a local asymptotic condition on the inverse of the transformation $z$ (see below \eqref{eq:c}).

We state first the result for $p_{a;x}$. 

\BT  \label{prop:AA}
Let  $(X,Y), \rho \in (-1,1)$ be a bivariate random vector with stochastic representation \eqref{eq:wzW}, where $z^*:[0,1] \to [0,1]$ is positive measurable function, $R\simeq F$, and let
$a\in (0,1]$ be a given constant. Suppose that there exists $\Ay \in (1, a/\abs{\rho} ) $ and for some $\ve\in (0,1)$ the function $z(x):=\rho x + z^*(x), x\in [0,1]$ is decreasing in $V_\ve:=[1/\Ay- \ve, 1/\Ay+ \ve]$ and $z(x)\le a/\Ay, \forall x\in  (1/\Ay, 1].$ Suppose that the inverse $z_\ve$ of $z$  in $V_\ve$ satisfies
\BQN\label{eq:c}
z_\ve(a/\Ay -  d/x)- 1/\Ay &\sim&  \frac{c d}{x}
 \EQN
locally uniformly for  $d>0$ with $c\in (0,\IF)$. If further $F\in GMDA(w)$ and \eqref{eq:cc} is satisfied with
$K_1:=-1/\Ay, K_2:=ca/\Ay$, then $\Ay$ is unique and
\BQN\label{eq:main:AA}
p_{a;x}&\sim&   \pk{I_1=I_2=1} \Gamma(\gamma_a+1)L_{K_1,K_2}(1/v(\Ayx))  \frac{\overline{F} (\Ayx)}{( v(\Ayx))^{\gamma_a} },
\EQN
where $\Ayx:= \Ay x, v(x):= xw( x), x>0$.
\ET
{\bf Remarks:} (a) If the random variable $W$ appearing in the stochastic representation \eqref{eq:wzW} possesses a positive density function $h$ continuous at
$1/\Ay$, then under the assumptions of \neprop{prop:AA} the asymptotics in \eqref{eq:cc} holds for any
$K_1 < K_2 $ with $\gamma_a=1$ and
\BQN \label{eq:gamma:LK12}
L_{K_1,K_2}(u)= (K_2 -K_1) h(1/\Ay), \quad u>0.
\EQN

(b) In view of \eqref{eq:resn} the tail asymptotics of $p_{a;x}$ given by \eqref{eq:main:AA} is faster than $\overline{F}(x)$. In fact for any constant $\mu >0$ we have
\BQNY
\limit{x}\frac{p_{a;x}}{(v(x))^\mu\overline{F}(x)}&=& 0.
\EQNY
Recall that the assumption $F \in GMDA(w)$ implies $\limit{x} x w(x)=\IF$.

c) As it can be seen from the proof of \neprop{prop:AA} the local behaviour of $z$ at $1/a_\rho$ is crucial. Since we assume that
$z$ is a decreasing function in $V_\ve$ the asymptotic of $z$ in $(1/a_\rho, 1/a_\rho+ \ve)$ is controlled by the asymptotic relation \eqref{eq:c}. Another possibility for $z$ is to assume that it is increasing in $(1/a_\rho - \ve)$ and decreasing in $(1/a_\rho - \ve)$ so that $1/a_\rho$ is a locally maximum for $z$. In this case we can still find the asymptotics of $p_{a;x}$, provided that additionally we assume that $z_\ve^{-}$ and $z_\ve^{+}$ are the inverses of $z$ in $(1/a_\rho- \ve, 1/a_\rho)$ and $(1/a_\rho, 1/a_\rho+ \ve)$, respectively such that
$$ z_\ve^{-}(a/a_\rho - 1/x) - 1/a_\rho \sim - \frac{c_{-}}{x},  \quad z_\ve^{+}(a/a_\rho - 1/x) - 1/a_\rho \sim \frac{c_{+}}{x}, \quad c_{-}, c_{+}\in (0,\IF)$$
locally uniformly in $x>0$.

In order to approximate the joint conditional excess distribution we need an asymptotic approximation as in the previous section of $\paxD$. In the setting of Model B we can approximate another quantity, namely $\paxDD$
defined by
$$ \paxDD:= \pk{ X> x [1+  \delta /v (a_\rho x)], Y> ax [1+  \eta /v (a_\rho x)]}, \quad x\to \IF$$
with $\delta, \eta\in [0,\IF)$ and $a_\rho$ as above. Note that the difference to $\paxD$ is that above we employ the normalisation function
$v(a_\rho x)$ and not $v(x)=x w(x)$. From the application point of view 
considering $\paxDD$ and not $\paxD$ is no restriction since
the interest is to be able to approximate the joint conditional excess function utilising some normalisation function. 
However, estimating $v(a_\rho x)$ leads to complication since also $\a_\rho$ need to be estimated.


\BT  \label{prop:AAA}
Under the assumptions and notation of \neprop{prop:AA} if $F$ has a density function $f$ continuous at $1/a_\rho$, then $\Ay$ is unique and
if further
\BQN\label{eq:cb}
z_\ve(a/\Ay +  d/x)- 1/\Ay &\sim&  - \frac{c d}{x}
 \EQN
locally uniformly for $d>0$, then for any $\eta, \delta\in [0, \IF)$ we have
\BQN\label{eq:main:AAAc}
\paxDD&\sim&   \pk{I_1=I_2=1} \frac{ h(1/a_\rho)}{a_\rho} (ca+1) \exp\Bigl( -\frac{ ca \eta + \delta}{ca+1}\Bigr)
\frac{\overline{F} (\Ayx)}{ v(\Ayx) },
\EQN
locally uniformly in $\delta,\eta$.
\ET

{\bf Remark:} If we drop the condition \eqref{eq:cb}, then the claim of \neprop{prop:AAA} still holds, provided that
$\delta\ge \eta \ge 0$.
\COM{can derive the asymptotics of
 $p_{a;\delta, \eta,x}$ also for the case $\eta > \delta \ge 0$, provided that additionally
\BQN\label{eq:cb}
z_\ve(a/\Ay +  1/x)- 1/\Ay &\sim&  - \frac{c}{x},
 \EQN
locally uniformly in $x>0$.  More specifically we have for any $\eta, \delta\in [0, \IF)$
\BQN\label{eq:main:AAAc}
p_{a;\delta, \eta, x}&\sim&   \pk{I_1=I_2=1} \frac{ h(1/a_\rho)}{a_\rho} (ca+1) \exp\Bigl( -\frac{ ca \eta + \delta}{ca+1}\Bigr)
\frac{\overline{F} (\Ayx)}{ v(\Ayx) },
\EQN
}

We present next two  examples.

{\bf Example 4.} Let $(X,Y)$ be a bivariate scale mixture random vector with stochastic representation \eqref{eq:wzW} where $\rho=0$.
Consider the function $z^*$ given by
$$z^*(x):= (1-   \abs{x}^p)^{1/p}, \quad p \in (0,\IF),x\in [-1,1].$$
The inverse function of $z:= z^*$ is $z^{-1}(y)= (1- y^p)^{1/p}, y\in [0,1]$. For any $a\in (0,1]$ the equation
$$ z^{-1}(a/s)= 1/s,  \quad s\in (1,\IF)$$
 has the unique solution $ \Ay:= (1+ a^p)^{1/p}\in (1, \IF)$. Furthermore \eqref{eq:c} and \eqref{eq:cb} hold with $c=a^{p-1}$.\\
Let  $W>0$ with positive density function $h$ being further independent of the positive random variable $R\simeq F$. If
 $F\in GMDA(w)$, then by \neprop{prop:AAA} 
\BQN\label{eq:impl}
\paxD&\sim&
 \pk{I_1=I_2=1}   \Ay^{p-2}  h(1/\Ay) \exp\Bigl( - \frac{\delta + a^p \eta}{1+a^p} \Bigr)  \frac{\overline{F} (\Ay  x)}{x w(\Ay x)}.
\EQN
Note that if $(RI_1W, I_2 R z^*(W))$ is a generalised symmetrised Dirichlet random vector, then $I_1,I_2, R,W$
are independent and $W$ possesses the density function $h(x)= p x^{p-1} g(x^p)$ with $g$ the density function of $beta(\alpha,\beta)$.

{\bf Example 5.} Under the setup of Example 4, with motivation from the dependence structure of elliptical random vectors we redefine $z^*$ as 
$$z^*(x):= \Orho \sqrt{1-   x^2}, \quad z(x):=\rho x+ z^*(x), \quad \rho,x\in (-1,1), \quad \rho_*:=\sqrt{1- \rho^2}.$$
First note that $z(\rho)=1$ is the maximal value of $z(x)$ for any $x\in [-1,1]$. Hence in order to apply \eqref{eq:main:AAAc} necessarily $a\in (\rho ,1]$. It can be easily checked that the assumptions of \neprop{prop:AA} are satisfied for
$\Ay:= \sqrt{1- 2a\rho+ a^2}/\rho_*$,
and \eqref{eq:c} holds with $c:=(a-\rho)/(1-a\rho)\in (0,\IF)$.
%
 Note further that $\Ay < a /\abs{\rho}$ and also \eqref{eq:cb} holds.
 In view of \eqref{eq:main:AAAc}  we obtain
\BQN\label{eq:1}
\paxDD&\sim&
\pk{I_1=I_2=1} \frac{\rho_*^2 h(1/\Ay)}{ 1- a \rho} \frac{\overline{F} (\Ay  x)}{x w(\Ay x)}
 \exp\Bigl( -\frac{ a^2 \eta +\delta- a\rho(\eta+ \delta) }{\Orho^2\Ay^2}\Bigr) .
\EQN
In the special case that $W^2 \simeq beta(1/2,1/2)$ and $\pk{I_i=1}=1/2, i=1,2$ with $I_1,I_2$ independent (the bivariate random vector $(X,Y)$ is elliptical distributed) we have
\BQNY
h(1/\Ay)&= &\frac{2 \Ay}{ \pi \sqrt{ \Ay^2- 1}}=  \frac{2 \Ay(1- \rho^2)}{a- \rho}.
\EQNY
Consequently   \eqref{eq:1} reduces to \eqref{eq:2} if additionally $\delta=\eta=0$.

\section{Three Applications}
Let $(X,Y)$ be a given bivariate random vector.  For some high threshold $x$ the approximation of the joint conditional excess random vector
$$ (X^{[x]}, Y^{[ax]}):= (X - x, Y -ax )\lvert X> x, Y> ax, \quad x \in (0, \IF), \quad a\in (0,1]$$
is of some interest in statistical applications if in particular suitable norming constants can be found so that the distribution of
$(X^{[x]}, Y^{[ax]})$ can be approximated by some known distribution function.\\
Another interesting problem of the bivariate extreme value theory is the asymptotic independence of $X$ and $Y$ if both have distribution functions in some max-domain of attractions. When $X$ and $Y$ are asymptotically independent an interesting 
topic also for application (see e.g., de Haan and Ferreira (2006)) is the estimation of the residual dependence index $\eta$. 
In our last application we give an explicit formula for $\eta$. 

In the light of our findings above we are able to discuss alternative solutions to both problems 
for the models of Section 2 and 3.

\subsection{Asymptotics of Conditional Excess Distribution}
 We start by considering the model of Section 2.  For any $s,t$ positive and some positive scaling function $w$ we have
$$ \pk{X^{[x]}> s/w(t), Y^{[ax]}> t/w(t)} = \frac{p_{a; s,t, x}}{p_{a; 0,0, x}}, \quad x>0.$$
If $x$ tends to infinity the asymptotics findings of Section 2 to approximate the above ratio.
More precisely, under the assumptions of \netheo{theo1} we have
\BQNY
\frac{p_{a; s,t, x}}{p_{a; 0,0, x}} \sim \frac{J_{s,t}}{J_{0,0}}.
\EQNY
By the definition $J_{s,t}$ depends on the limit function $\xi_a$. Denote by $(E_1,E_2)$ a bivariate random vector with positive components and
survivor  function given by $J_{s,t}/ J_{0,0},s,t\in (0,\IF)$.
Then the above asymptotics can be cast into joint convergence in distributions, namely
if $(X_n, Y_{n}), n\ge 1$ is a sequence of bivariate random vectors defined in the same probability space such that $(X_n,Y_n)\equaldis (X^{[n]},Y^{[an]}), n\ge 1$, then  we have the convergence in distributions
\BQN \label{eq:cup}
( g(n)  X^{[n]}, g(n) Y^{[an]}) &\todis & (E_1,E_2), \quad n \to \IF,
\EQN
where the scaling function $g$ equals $w$.\\
The limiting random vector has distribution function which clearly depends on $a$. Further, $E_1$ and $E_2$ can be dependent for instance in the setup of Example 3 taking $a=1$. In the next model this joint distribution of $(E_1,E_2)$ is a product distribution
which seems to be more relevant for statistical applications. \\
 Assume next that $(X,Y),a, a_\rho,\rho$ satisfy the assumptions of \neprop{prop:AAA}.
 For any $s,t$ positive \eqref{eq:main:AAAc} implies  (set $\bar s:= \Ay s, \tilde t := \Ay t/a$)
\BQNY
\lefteqn{ \frac{ \pk{ X> x + s/w(\Ay x), Y> ax+ t/w(\Ay x)}} { \pk{ X> x, Y> ax }} }\\
&=&
\frac{ \pk{ X> x(1+ \bar s/v(\Ay x)), Y> ax(1+ \tilde t/v(\Ay x))}}  {\pk{X> x, Y> a x}}\\
&=&
\frac{ p_{a, \bar s, \tilde t, \rho; x} }{p_{a;x}}\\
&\sim & \exp\Bigl( -\frac{ \Ay s+ ca t \Ay /a  }{ca+1}\Bigr)
=: \exp(-  sD_{a,c}- t D_{a,c}^*),
\EQNY
where
$$D_{a,c}:=\frac{\Ay}{ca+1}, \quad  D_{a,c}^*:= \frac{ \Ay c}{a(ca+1)}.$$
%
\COM{
If additionally this model we cannot in general determine the tail asymptotic behaviour of $X$ and $Y$ since we do not make explicit assumptions on the upper tail asymptotics of $U_1$ and $U_2$. On the other side, in theatdo not
 and assume that \eqref{eq:gamma:LK12} holds. For any $x$ large we may write ($\bar x:= x+ t/[a w(\Ay x)]$)
\BQNY
\lefteqn{\pk{ X> x + s/w(\Ay x), Y> ax+ t/w(\Ay x)}}\\
&=&
\pk{ X>  \bar x + [s- t/a]/w(\Ay x), Y> a \bar x }, \quad \forall s,t\inr.
\EQNY
Next, applying \eqref{eq:main:AAB} with $\delta= s- t/a$ we obtain
\BQNY
\frac{\pk{ X> x + s/w(\Ay x), Y> ax+ t/w(\Ay x)}}
{\pk{ X> x, Y> ax}} &\sim  & \frac{ \overline{F} (\Ay \bar x)}
{ \overline{F} (\Ay x)} \exp(- [s- t/a]\Ay/(ac+1))\\
&=& \exp(- \Ay t/a) \exp(- [s- t/a]\Ay/(ac+1))\\
&=& \exp(-  s\Ay/(ac+1)- t \Ay[1-1/(ac+1)]/a)\\
&=:& \exp(-  sD_{a,c}- t D_{a,c}^*),
\EQNY
where
$$D_{a,c}:=\Ay/(ac+1), \quad  D_{a,c}^*:= \Ay[1-1/(ac+1)]/a.$$}

Consequently, with $(X_n, Y_{n}), n\ge 1$ as defined above \eqref{eq:cup}
holds with
$$ g(x)= w(\Ay x), \quad x> 0$$
and $E_1,E_2$ two independent exponential random variables with mean $1/D_{a,c}$ and $1/D_{a, c}^*$, respectively.

Under the setup of Example 5
$$ D_{a,c}:= \frac{1- a \rho}{ \Ay(1- \rho^2)}, \quad  D_{a,c}^*:=  \frac{a -\rho}{ \Ay(1- \rho^2)}.   $$
Thus the convergence in distribution in \eqref{eq:cup} holds in particular if
$(U_1,U_2)$ is uniformly distributed on the unit circle of $\R^2$.
We note that the approximation of the conditional excess distribution we do not assume a specific tail asymptotics of $X$ and $Y$.

\subsection{Asymptotic Independence and Max-domain of Attraction}
A common measure of the asymptotic dependence between $(X,Y)$ is
the tail dependence function
$$\limit{x} \frac{\pk{G_1(X)> 1- s/x, G_2(Y)> 1- t/x}}{\min( \pk{G_1(X)> 1- s/x}, \pk{G_2(Y)> 1- t/x})}:=l(s,t), \quad s,t\in (0,\IF) $$
(when it exists) where $G_1,G_2$ are the distribution functions of $X$ and $Y$, respectively.   If $l(1,1)=0$, then we say that $X$ and $Y$ are asymptotically independent. See for instance de Haan and Ferreira (2006), Reiss and Thomas (2007), H\"usler and Li (2009),
 Das and Resnick (2009), or Peng (2010) for more details concerning modelling of asymptotic independence in the context of extreme values.\\
We discuss briefly the asymptotic independence for scale mixture distributions with $(U_1,U_2)$ specified by our Model A.
It can be seen by Example 1 that for particular $U_1,U_2$ the limit $l(s,t)$ can be positive, thus asymptotic independence does not hold.
However, under the setup of Example 2 \eqref{FGM} implies that $l(s,t)=0,\forall s,t\in(0,\IF)$, and thus $X$ and $Y$ are asymptotically independent and both
$X$ and $Y$ have distribution function in the Gumbel max-domain of attraction.

We deal next with Model B assuming that $(X,Y)$ has stochastic representation \eqref{eq:wzW} with $\rho\in [0,1)$. The case $\rho \in (-1,0)$ follows with similar arguments.\\
In the following we specify the asymptotic behaviour of $W$ and $z(W)$. Explicitly, we assume that for some $\gamma_1,\gamma_2\in [0,\IF)$
$$ \limit{x}\frac{\pk{W> 1- s/x}}{\pk{W> 1- 1/x}}= s^{\gamma_1}, \quad  \limit{x} \frac{\pk{z(W)> 1- s/x}}{\pk{z(W)> 1- 1/x}}= s^{\gamma_2}, \quad \forall s> 0.$$
As in Example 1 applying \eqref{eq:berm} we obtain
$$ \pk{R W> x}\sim  \Gamma(\gamma_1+1) \pk{W> 1 - 1/\vx}\overline{F}(x)$$
and
$$ \pk{R z(W)> x}\sim  \Gamma(\gamma_2+1) \pk{z(W)> 1 - 1/\vx}\overline{F}(x).$$
Next set $z(x):= \rho x+ z^*(x) \le  1, \forall x\in [-1,1] $
and assume that  $z^*(x) \le b < 1, \forall x \in [0,1]$. Applying \eqref{eq:resn} we obtain
\BQNY
\pk{Y> x}&=& \pk{I_1=1,I_2=-1} \pk{R(\rho  W - z^*(W))> x}\\
&&+
\pk{I_1=-1,I_2=1} \pk{R(-\rho  W + z^*(W))> x}\\
&&+
 \pk{I_1=1, I_2=1} \pk{R(\rho  W + z^*(W))> x}\\
&\sim&\pk{I_1=1,I_2=1} \pk{R z(W)> x}\\
&\sim&\pk{I_1=1, I_2=1} \Gamma(\gamma_2+1) \pk{z(W)> 1 - 1/\vx}\overline{F}(x).
\EQNY
Similarly,
\BQNY
 \pk{X> x}&\sim&   \pk{I_1=1} \Gamma(\gamma_1+1) \pk{W> 1- 1/\vx}\overline{F}(x).
 \EQNY
Consequently, in view of \eqref{eq:wb} both $X$ and $Y$ have distribution functions in the Gumbel max-domain of attraction with the same scaling function $w$. Let $b_i(x),i=1,2$ be defined asymptotically by
$$ b_i(x):= G_i^{-1} (1- 1/x), \quad x> 1,$$
where $G_i^{-1}$ is the generalised inverse of $G_i,i=1,2$. In view of \eqref{eq:resn} we have
\BQN\label{eq:resn2}
\limit{x} \frac{b_1(x)}{b_2(x)}&=& 1.
\EQN
Furthermore (see e.g., Falk et al.\ (2004))
\BQN\label{eq:resn3}
w(b_i(x)) [G_i^{-1}(1- s/x)- b_i(x) ] &=& - \ln s, \quad \forall s\in (0,\IF).
\EQN
If $X,Y$ are such that the conditions of \neprop{prop:AA} hold, then comparing the asymptotics of $\pk{X> b_1(x/s), Y> b_2(x/t)}$ and $\pk{X> b_1(x/s)},\pk{Y> b_2(x/t)}$ we obtain utilising further \eqref{eq:resn} and \eqref{eq:resn2}
$$ l(s,t)=0, \quad \forall s,t\in (0,\IF).$$
Consequently, $X$ and $Y$ are asymptotically independent with distribution function in the max-domain of attraction of a bivariate distribution with unit Gumbel marginals which is a product distribution.

\subsection{Residual Tail Dependence}

Modeling of dependence and asymptotic dependence is often done in the framework of copula, where the marginal
distributions are transformed. The asymptotic dependence does not change under monotone transformation of
marginal distributions. For $X,Y$ with asymptotically independent components it is of some interest to quantify the asymptotic independence in terms of some measures.
Let $G_1,G_2$ be the distribution function of $X$ and $Y$, respectively.
One successful approach to model the asymptotic independence is the estimation of the residual dependence index
$\eta \in (0,1]$  (see  Peng (1998,2008,2010), de Haan and Ferreira (2006), Hashorva (2010)). So if for some  $x,y$ positive
$$ S_u(x,y):= \frac{\pk{G_1(X)> 1- x/u, G_2(Y)> 1- y/u}}{\pk{G_1(X)> 1- u, G_2(Y)> 1- u}}\to
 S(x,y), \quad u \to \IF,$$
then for any $c>0$
$$ S(cx,cy)= c^{1/\eta} S(x,y)$$
and the function $ S_u(1,1)$  is regularly varying with index $-1/\eta$.
Other authors refer to $\eta$ as the coefficient of tail dependence (Ledford and Tawn (1998), Resnick (2007), Reiss and Thomas (2007)).

As mentioned above in Model A asymptotic independence is not always observed as for instance in the setup of Example 1.
However, as noted above for Example 2 asymptotic independence is observed. We calculate $\eta$ for that  example.
Denote next by $G^{-1}_i,i=1,2$ the generalised inverse of $G_i,i=1,2$. 
Since further $\limit{u} b_1(u)=\limit{u} b_2(u)= \IF$, by \eqref{eq:resn3} we can write for any $x,y \in (0,1)$ as $u\to \IF$
\BQN\label{merzi}
\limit{u} \frac{ S_u(x,y)}{S_u(1,1)}&=& \limit{u} \frac{\pk{G_1(X)> 1- x/u, G_2(Y)> 1- y/u}}{\pk{G_1(X)> 1- 1/u, G_2(Y)> 1- 1/u}}\notag \\
&=& \limit{u} \frac{\pk{X> G^{-1}(1- x/u), Y> G_2^{-1}(1- y/u)}}{\pk{X> G_1^{-1}(1- 1/u), Y > G_2^{-1}(1- 1/u)}}\notag \\
&=& \limit{u} \frac{\pk{X> b_1(u) - \ln x / w(b_1(u)), Y> b_2(u) - \ln y / w(b_2(u))}}{\pk{X> b_1(u) , Y > b_2(u)}}\notag \\
&=& \limit{u} \frac{\pk{X> b_1(u) - \ln x / w(b_1(u)), Y> b_1(u) - \ln y (1+o(1))/ w(b_1(u))}}{\pk{X> b_1(u) , Y > b_1(u)+ o(1)/w(b_1(u))}}.
\EQN
Hence by \eqref{FGM} we obtain
\BQN
\limit{u} \frac{ S_u(x,y)}{S_u(1,1)}&=:&S(x,y)=\frac{1}{\Gamma(\gamma_1+\gamma_2+1)}\int_0^{\IF}(t+ \ln x)_+^{\gamma_1}(t+ \ln y)_+^{\gamma_2}\exp(-t)\, dt.
\EQN
Consequently, since for any $c>0$ we have $S(cx,cy)=c S(x,y)$ we conclude that
$$ \eta= 1.$$

We consider next Model B. Let therefore $(X,Y)$ be as in our second application
satisfying further the assumptions of \neprop{prop:AAA}. Since  $X,Y$ are asymptotically independent we deal next with the
calculation of residual dependence index $\eta$. We assume  the scaling function $w$ (see the assumptions of our second application) is such that
\BQN\label{eq:wro}
\limit{u} \frac{ w(cu)}{w(u)}&=& c^{\lambda-1}, \quad \forall c\in (0,\IF).
\EQN
Since necessarily $\limit{u} u w(u)= \IF$ we require further that $\lambda\in [0,\IF)$. Thus $w(x)= x^{\lambda-1} L(x)$ with $L$ a positive slowly varying function at infinity. If $\lambda=0$ we assume further that $\limit{u} L(u)= \IF$.

\underline{\bf Case $\limsup_{u\to \IF} w(u)< \IF$}:\\
In view of \eqref{eq:resn2} we have for any $y\in (0,\IF)$
$$b_2(u) - \ln y /w(b_2(u))= b_1(u) + o(1)- (1+o(1))\ln y / w(b_1(u)) = b_1(u)- (1+o(1)) \ln  y / w(b_1(u)), \quad u\to \IF.   $$
Since further $\limit{u} b_1(u)= \IF$, \eqref{merzi} and \eqref{eq:wro} imply for any $x,y \in (0,1)$
\BQNY
\limit{u} \frac{ S_u(x,y)}{S_u(1,1)}
&=& \limit{u} \frac{\pk{X> b_1(u) -  \alpha_\rho^{\lambda- 1}\ln x / w(\alpha_\rho b_1(u)), Y> b_1(u) - \alpha_\rho^{\lambda- 1} \ln y (1+o(1))/ w(\alpha_\rho b_1(u))}}{\pk{X> b_1(u) , Y > b_1(u)+ o(1)/w(\alpha_\rho b_1(u))}}.
\EQNY
As in our first application we obtain
\BQNY
\limit{u} \frac{ S_u(x,y)}{S_u(1,1)}&=& \exp( \alpha_\rho ^{\lambda- 1} [D_{a,c} \ln x+ D_{a,c}^* \ln y]).
\EQNY
Consequently
$$\eta^{-1}=   \alpha_\rho ^{\lambda- 1} [D_{1,c} + D_{1,c}^* ]=   \alpha_\rho ^{\lambda} \frac{1+ c}{(c+ 1)}=\alpha_\rho ^{\lambda} .$$
Since $\alpha_\rho > 1,$ then clearly $\eta \in (0,1]$. It is interesting that $\eta$ depends only on $\alpha_\rho$ and $\lambda$ and not on $c$.

\underline{\bf Case $\limit{u} w(u) = \IF$}:\\
In order to calculate $\eta$ we need to assume further a certain relationship between $b_1(u)$ and $b_2(u)$. In view of \eqref{eq:resn2} and the assumption on $w$
suppose further that $b_2$ is such that
\BQN
\limit{u} w(b_2(u))[ b_2(u)- b_1(u)]&=& \xi\inr .
\EQN
As above for any $x,y\in (0,\IF)$ we obtain as $u\to \IF$
\BQNY
\limit{u} \frac{ S_u(x,y)}{S_u(1,1)}&=& \limit{u} \frac{\pk{X> b_1(u) -  \alpha_\rho^{\lambda- 1}\ln x / w(\alpha_\rho b_1(u)),
Y> b_1(u) -  (1+o(1))\alpha_\rho^{\lambda- 1} (\ln y+ \xi)/ w(\alpha_\rho b_1(u))}}{\pk{X> b_1(u) , Y > b_1(u)+ (1+o(1)) \alpha_\rho ^{\lambda-1}\xi /w(\alpha_\rho b_1(u))}}\\
&=& \exp( \alpha_\rho ^{\lambda- 1} [D_{a,c} \ln x+ D_{a,c}^* \ln y]),
\EQNY
hence again $\eta= \alpha_\rho^{- \lambda}\in (0,1]$. When $\lambda=0$, then $\eta=1$, otherwise we have $\eta \in (0,1)$. Note in passing that neither $\xi$ nor $c$ appear in the expression of the residual tail dependence $\eta$.

For statistical models estimation of $\eta$ is important. In view of our derivation for this model we can estimate $\eta$ by estimating first $\alpha_\rho$ and then $\lambda$. An estimation of $\lambda$ can be obtained as in Hashorva (2010), whereas estimation of $\alpha_\rho$ is not as straightforward. In the more specific model of Example 5 $\alpha_\rho$ can be estimated if we estimate $\rho$. Estimation of $\alpha_\rho$ will be discussed in a forthcoming paper.

\section{Proofs}
\proofprop{theo1}  Since $R$ is independent of the bivariate random vector $(U_1, U_2)$ and $U_1 \le 1$ almost surely for any
$\delta,\eta \in [0,\IF)$ we have
\BQNY
\paxD&= &\pk{ R U_1> x(1+ \delta/\vx), R U_2 > a x(1+ \eta/\vx)}\\
& =& \int_x^{\IF}  \pk{  U_1> x/r,  U_2 > a x/r}\, dF(r), \quad \forall x>0.
 \EQNY
Let $\ve$ be a positive constant. The assumption \eqref{eq:paxg} implies that for any constant $c \in (1, 1+\ve)$ we have $\pk{U_1> 1/c, U_2> a/c}\in (0,1)$. In view of \eqref{eq:resn}
 \BQN
 \label{eq:rap}
 \limit{x} \frac{\overline{F}(\alpha x)}{\overline{F}(x)}&=&0, \quad \forall \alpha \in (1,\IF)
 \EQN
 and the fact that $(U_1, U_2)$ is independent of $R$ we obtain
\BQN\label{eq:eps}
\paxD& \sim& \int_{ x(1+ \delta/\vx)} ^{cx} \pk{  U_1>x(1+ \delta/\vx)/r, U_2> a x(1+ \eta/\vx)/r }\, dF(r)\\
&=& \int_{\delta} ^{ \vx (c-1) } f(s,\vx, \delta,\eta)\, dF(s/w(x)+ x),
 \EQN
 with
 $$ f(s,\vx, \delta,\eta):= \pb{  U_1>\frac{1+ \delta/\vx}{1+ s/\vx}, U_2> \frac{a (1+ \eta/\vx)}{1+ s/\vx} }, \quad s\ge \delta.$$
If  $\delta=\eta=0$, then \eqref{eq:paxg} implies that
$$f(s,1, 0,0)=\pk{  U_1>1/(1+ s), U_2> a /(1+ s) }, \quad s>0   $$
is regularly varying at 0. By the max-domain of attraction assumption on $F$ and \eqref{eq:uv} the result for $\delta=\eta=0$
follows easily applying further Potters bounds (see e.g., de Haan and Ferreira (2006)) for the integrand and utilising Lemma 7.5 and 7.7 in Hashorva (2007). The general case $\delta$ or $\eta$ can be established utilising the result for $\delta=\eta=0$, \eqref{theo1:1:1} and the fact that
 $f(s,\vx, \delta,\eta)\le f(s,\vx, 0,0), \delta,\eta\ge0$, and thus the proof is complete. \QED
\def\vyyx{v(\Ayx)}

\proofprop{prop:AA} Define next
$$x_*:=\Ay x, \quad v(x)= xw(x), \quad F_x(s):= F(x^*[1+ s/\vyyx]), \quad s\inr, x>0.$$
By the independence of $I_1,I_2$ and $RW$ we may write for any $x>0$
\BQNY
p_{a;x}&=&\pk{RW> x,\rho RW  - z^* (W)> ax} \pk{I_1=1, I_2=-1}\\
&& +
\pk{RW> x,R z(W)> ax} \pk{I_1=1,I_2=1}\\
&=:&J_1 \pk{I_1=1, I_2=-1}+ J_2\pk{I_1=1,I_2=1}.
\EQNY
If $\rho \le 0$, then the fact that $z^*$ is non-negative  implies $p_{a;x}=J_2.$
When $\rho \in (0,1)$ the assumptions $a/\rho > \Ay$ yields further
\BQNY
J_1 \le \pk{RW > a/\rho x} \le \pk{R> a/\rho x}= \overline{F}(a/\rho x).
\EQNY
Since $z(s)\le  a/\Ay$ for any $s\in [1/\Ay ,1]$ we have
\BQNY
\pk{W> x/r, z(W)> ax/r}& \le &
\pk{W>1/\Ay , z(W)> ax/r}\\
&\le &\pk{W> 1/\Ay, z(W)> a/\Ay}\\
&=&0, \quad \forall x,r>0, x\le  r\le  \Ay x.
\EQNY
Hence for any  $\ve\in (0,1)$
\BQNY
J_2&=&\int_{\Ay x }^{(\Ay+\ve) x} \pk{W> x/r,z(W)> ax/r}\, dF(r)+ \int_{(\Ay+\ve) x}^{\IF} \pk{W> x/r,z(W)> ax/r }\, dF(r).
\EQNY
Since (by the assumption) the function $z$ is decreasing and possesses an inverse function $z_\ve$ in $[1/\Ay- \ve, 1/\Ay+\ve]$ for some
given $\ve>0$,  we have  
\BQNY
\lefteqn{\int_{\Ay x }^{(\Ay+\ve) x} \pk{W> x/r,z(W)> ax/r}\, dF(r)}\\
&=&\int_{0}^{\ve \vyyx} \pk{W> 1/\Ay [1 + s/\vyyx]^{-1} ,z(W)> a/\Ay [1 + s/\vyyx]^{-1} }\, dF_x(s)\\
&=&\int_{0}^{\ve \vyyx}\pk{W> 1/\Ay [1 + s/\vyyx]^{-1} ,W < z_\ve( a/\Ay [1 + s/\vyyx]^{-1}) }\, d F_x(s)\\
&=&\int_{0}^{\ve \vyyx} \pk{ c a s/(\Ay \vyyx)(1+o(1))     >  W- 1/\Ay >  - s/(\Ay \vyyx)(1+o(1))   }\, d F_x(s).
\EQNY
Hence by the assumptions on $W$ and $F$ applying Potters bound for the integrand and utilising Lemma 7.5 and 7.7 in Hashorva (2007)  we obtain
\BQNY
\lefteqn{\int_{\Ay x }^{(\Ay+\ve) x} \pk{W> x/r,z(W)> ax/r}\, dF(r)}\\
&\sim& \int_{0}^{\IF} s^{\gamma_a}\exp(-s)\, ds
L_{-1/\Ay,c a/\Ay   }(1/\vyyx) \frac{ \overline{F} (x_*)}{( \vyyx)^{\gamma_a} } .
\EQNY
In view of \eqref{eq:rap} $\Ay$ is necessarily unique, hence applying \eqref{eq:resn} as $x\to \IF$
\BQNY
p_{a;x}&= & \overline{F}(a/\rho)+ (1+o(1))\pk{I_1=I_2=1}\Gamma(\gamma_a+1)
L_{-1/\Ay,c a/\Ay   }(1/\vyyx) \frac{ \overline{F} (x_*)}{( \vyyx)^{\gamma_a} }
+ O( \overline{F}((\Ay+ \ve)x) \\
&\sim&\pk{I_1=I_2=1}\Gamma(\gamma_a+1) L_{-1/\Ay,c a/\Ay   }(1/\vyyx)  \frac{ \overline{F} (x_*)}{( \vyyx)^{\gamma_a} },
\EQNY
and thus the result follows.  \QED

\proofprop{prop:AAA} By the assumption on the density function $h$ of $W$ we have that \eqref{eq:cc} holds for any $K_1< K_2, K_1,K_2\inr$ with $\gamma_a=1$. As in the proof above for any $\delta \in [0,\IF)$ and $\ve>0$ small enough  we obtain (set $\xi:=(ca \eta +\delta)/(ca+1)$)
\BQNY
\paxDD
&=&(1+o(1))\pk{I_1=I_2=1}\\
&& \times \int_{\xi }^{\ve \vyyx} \pk{ c a (s- \eta)/(\Ay \vyyx)(1+o(1))     >  W- 1/\Ay >  (\delta- s)/(\Ay \vyyx)(1+o(1))   }\, d F_x(s)\\
&\sim & \pk{I_1=I_2=1} \int_{\xi}^\IF [(ca +1)s - \delta -\eta] \exp(-s) \, ds
\frac{ h(1/a_\rho)}{a_\rho} \frac{ \overline{F} (x_*) }{ \vyyx }, \quad x\to \IF,
\EQNY
hence the result follows.  \QED

We conclude this section with the proof of \eqref{eq:apend:1}.\\
 For all $x$ large and $s> 0$ by the independence of $S_1$ and $S_2$ for any $\delta_i\in [0,\IF),i=1,2,s>0$ we may write
(set $G_{2,x}(z):=G_2(1- z/x),  s_i:=s- \delta_i, s,x,z\in (0,\IF)$)
\BQNY
\lefteqn{ \pk{U_1 > 1- s_1/x, U_2 > 1- s_2/x} }\\
&=& \int_0^1 \pk{ \lambda_i S_1 > 1- s_i/x - \overline{\lambda_i} y, i=1,2}\, d G_2(y)\\
 &=& \int_0^{\IF} \pk{ \lambda_i S_1 > 1- s_i/x - \overline{\lambda_i} (1- z/x), i=1,2}\, d G_{2,x}(z)\\
  &=& \int_0^{\IF} \pk{ S_1 > 1-  [ s_i- \overline{\lambda_i} z]/(x \lambda_i) , i=1,2}\, d G_{2,x}(z)\\
  &=&  \prod_{i=1}^2 \overline{G}_i(1- 1/x)\int_0^{\IF} \frac{\pk{ S_1 > 1-  [ s_i- \overline{\lambda_i} z]/(x \lambda_i) , i=1,2}}
  {\overline{G}_1(1- 1/x)}\, d G_{2,x}(z)/\overline{G}_2(1- 1/x).
  \EQNY
The asymptotic behaviour of $\overline{G}_i,i=1,2$ implies
\BQNY
\lefteqn{ \pk{U_1 > 1- s_1/x, U_2 > 1- s_2/x} }\\
&\sim &\prod_{i=1}^2 G_i(1- 1/x)\gamma_2  \int_0^\IF
\Bigl( \max \Bigl(0,   \min (   [ s_1- \overline{\lambda_1}z]/\lambda_1, [ s_2- \overline{\lambda_2}z]/ \lambda_2)\Bigr)
\Bigr)^{\gamma_1} z^{\gamma_2- 1} \, dz.
  \EQNY
Now, for $\delta=\eta=0$ we may write further
\BQNY
\lefteqn{ \pk{U_1 > 1- s/x, U_2 > 1- s/x} }\\
&\sim &\prod_{i=1}^2 G_i(1- 1/x)\gamma_2  \int_0^{\min (s/ \overline{\lambda_1}, s/\overline{\lambda_2})}
\Bigl( \min \Bigl[   [ s- \overline{\lambda_1}z]/\lambda_1, [ s- \overline{\lambda_2}z]/ \lambda_2\Bigr] \Bigr)^{\gamma_1}z^{\gamma_2- 1} \, dz.
  \EQNY
Since $\lambda_2 \ge \lambda_1>0$,  for any $0 \le z \le s$
$$ \frac{s}{\overline{\lambda_2}} \ge \frac{s}{\overline{\lambda_1}}, \quad
\frac{ s- \overline{\lambda_1}z}{\lambda_1 }\ge
\frac{ s- \overline{\lambda_2}z}{ \lambda_2}$$
hence as $x\to \IF$
\BQNY
\lefteqn{ \pk{U_1 > 1- s/x, U_2 > 1- s/x} }\\
&\sim& \prod_{i=1}^2 \overline{G}_i(1- 1/x)\gamma_2  \Biggl[\lambda_2^{- \gamma_1}  \int_0^{s}
( s- \overline{\lambda_2}z)^{\gamma_1} z^{\gamma_2- 1} \, dz
+ \lambda_1^{- \gamma_1} \int_s^{s/\overline{\lambda_1}}
(  s- \overline{\lambda_1}z)^{\gamma_1}z^{\gamma_2- 1} \, dz\Biggr].
\EQNY

\end{document}